\documentclass[12pt]{amsart}

\usepackage{amssymb}

\usepackage{enumitem}

\makeatletter
\@namedef{subjclassname@2020}{%
  \textup{2020} Mathematics Subject Classification}
\makeatother

\usepackage{url}

\newtheorem{theorem}{Theorem}[section]
\newtheorem{corollary}[theorem]{Corollary}
\newtheorem{lemma}[theorem]{Lemma}
\newtheorem{proposition}[theorem]{Proposition}

\theoremstyle{definition}
\newtheorem{definition}[theorem]{Definition}
\newtheorem{remark}[theorem]{Remark}
\newtheorem{example}[theorem]{Example}

\numberwithin{equation}{section}

\newcommand{\abs}[1]{\left|{#1}\right|}
\newcommand{\Z}{\mathbb Z}
\newcommand{\Q}{\mathbb Q}
\newcommand{\R}{\mathbb R}
\newcommand{\CC}{\mathcal C}
\newcommand{\NN}{\mathcal N}
\newcommand{\Gal}{\text{Gal}}
\newcommand{\ord}{\text{ord}}
\newcommand{\ds}{\displaystyle}

\begin{document}


\baselineskip=17pt


\title[Monogenic Cyclic Polynomials in Recurrence Sequences]{Monogenic Cyclic Polynomials\\ in Recurrence Sequences}

\author[J. Harrington]{Joshua Harrington}
\address{Department of Mathematics, Cedar Crest College, Allentown, Pennsylvania, USA}
\email{Joshua.Harrington@cedarcrest.edu}

\author[L. Jones]{Lenny Jones}
\address{Professor Emeritus, Department of Mathematics, Shippensburg University, Shippensburg, Pennsylvania 17257, USA}
\email{doctorlennyjones@gmail.com}


\date{}

\begin{abstract}
 Let $f(x)\in {\mathbb Z}[x]$ be an $N$th degree polynomial that is monic and irreducible over ${\mathbb Q}$. We say that $f(x)$ is {\em monogenic} if $\{1,\theta,\theta^2,\ldots ,\theta^{N-1}\}$ is a basis for the ring of integers of ${\mathbb Q}(\theta)$, where $f(\theta)=0$. We say that $f(x)$ is {\em cyclic} if the Galois group of $f(x)$ over ${\mathbb Q}$ is the cyclic group of order $N$. In this article, we investigate the appearance of monogenic cyclic polynomials in certain polynomial recurrence sequences.
\end{abstract}


\maketitle

\section{Introduction}\label{Section:Intro}
 Let $f(x)\in {\mathbb Z}[x]$ be an $N$th degree polynomial that is monic and irreducible over ${\mathbb Q}$. We define $f(x)$ to be {\em monogenic} if $\{1,\theta,\theta^2,\ldots ,\theta^{N-1}\}$ is a {\em power basis} for the ring of integers $\Z_K$, of $K=\Q(\theta)$, where $f(\theta)=0$. We say that   $f(x)$ is {\em cyclic} if $\Gal(f)\simeq C_N$, where $\Gal(f)$ denotes the Galois group of $f(x)$ over $\Q$, and $C_N$ denotes the cyclic group of order $N$. We also say that $K$ is a {\em degree-$N$ cyclic extension of $\Q$} in this case. We define two degree-$N$ monogenic cyclic polynomials $f(x)$ and $g(x)$, with $f(\alpha)=g(\beta)=0$, to be {\em equivalent} if $\Q(\alpha)=\Q(\beta)$, and {\em distinct} otherwise. We illustrate these concepts with the following example.
 \begin{example}\label{Ex}
   Let 
   \begin{align*}
     f(x)&=x^4-4x^2+2,\\
     g(x)&=x^4-8x^3+20x^2-16x+2\\
     & \qquad \mbox{and}\\
     h(x)&=x^4+4x^2+2.
   \end{align*} Note that $f(x)$, $g(x)$ and $h(x)$ are all 2-Eisenstein, and hence irreducible over $\Q$. Let 
   \[\theta_{j,k}:=(-1)^j\sqrt{2+(-1)^k\sqrt{2}}, \quad \mbox{where} \ j,k\in \{0,1\}.\]  Observe that $\theta_{j,k}\in \R$.   Using a computer algebra system, such as Magma, Maple or Sage, it is straightforward to verify the following information:
   \begin{itemize}
     \item the zeros of $f(x)$ are $\theta_{j,k}$,
     \item the zeros of $g(x)$ are $2+\theta_{j,k}$,
     \item the zeros of $h(x)$ are $i\theta_{j,k}$, where $i^2=-1$,
     \item $\Gal(f)\simeq \Gal(g)\simeq \Gal(h)\simeq C_4$,
     \item $\Delta(f)=\Delta(g)=\Delta(h)=2^{11}$,
     \item $f(x)$, $g(x)$ and $h(x)$ are all monogenic.
   \end{itemize}
   Then, it is easy to see that: 
   \begin{itemize}
     \item $f(x)$ and $g(x)$ are equivalent,
     \item $f(x)$ and $h(x)$ are distinct,
     \item $g(x)$ and $h(x)$ are distinct.
   \end{itemize}
 \end{example}

 In a series of papers in the mid-1980's, Marie-Nicole Gras proved several results concerning the existence of monogenic cyclic extensions of $\Q$. Among these results  was the following theorem.
\begin{theorem}{\rm \cite{Gras2}}\label{Thm:Gras}
  Let $\ell$ be a prime, and let $K$ be a degree-$\ell$ cyclic extension of $\Q$. If $\ell\ge 5$, then
$\Z_K$ does not have a power basis except in the case when $2\ell+1$ is prime and $K=\Q(\zeta+\zeta^{-1})$, the maximal real subfield of the cyclotomic field $\Q(\zeta)$, where $\zeta$ is a primitive $2\ell+1$ root of unity. 
\end{theorem}
To illustrate Theorem \ref{Thm:Gras} in the context of monogenic cyclic polynomials, suppose first that $\ell=7$. Then, Theorem \ref{Thm:Gras} tells us that no degree-7 monogenic polynomial $f(x)$ exists with $\Gal(f)\simeq C_7$, since $2\ell+1=15$ is not prime. However, if $\ell=5$, then $2\ell+1=11$ is prime, and we conclude from Theorem \ref{Thm:Gras} that there exists a degree-5 monogenic polynomial $f(x)$ with $\Gal(f)\simeq C_5$. Indeed, if $f(x)=x^5+x^4-4x^3-3x^2+3x+1$, which is the minimal polynomial over $\Q$ of $\zeta+\zeta^{-1}$, where $\zeta$ is a primitive $11$th root of unity, then it is easy to verify that $f(x)$ is monogenic and $\Gal(f)\simeq  C_5$. Moreover, Theorem \ref{Thm:Gras} says that $f(x)$ is essentially the \emph{only} such polynomial, in the sense that if $g(x)$ is a degree-5 monogenic polynomial with $g(x)\ne f(x)$ and $\Gal(g)\simeq C_5$, then $g(x)$ is equivalent to $f(x)$.  For example, it is easy to see that $g(x)=x^5-11x^4+44x^3-77x^2+55x-11$ is equivalent to $f(x)$ since $g(x)=-f(2-x)$.

In \cite{Gras3}, Gras extended the ideas in Theorem \ref{Thm:Gras} to the situation when $2\ell^r+1$ is prime with $r\ge 1$. In \cite{Gras1}, Gras proved another result that implies that there exist at most finitely many monogenic cyclic extensions of $\Q$ of degree $N$ when $\gcd(N,6)=1$. Because of this seeming rarity of monogenic cyclic extensions of $\Q$, one conclusion that could be drawn from the work of Gras is that we might expect some difficulty in determining monogenic cyclic polynomials defining these cyclic extensions. However, in some sense, that conclusion is incorrect since Gras gives us an indication of where to look for monogenic cyclic polynomials, namely maximal real subfields of cyclotomic extensions. On the other hand, the complete story is more complicated, especially when $\gcd(N,6)>1$. For example, it follows easily from \cite{JonesShanks} that there are infinitely many distinct monogenic cyclic polynomials of degree $N=3$. Several authors \cite{HJAA,HJBAMS,JonesQuarticsBAMS,JonesAA,JonesRecipQuartics,MNSU,Voutier} have contributed to an examination of the case $N=4$. Most recently, it has been shown in \cite{Voutier}, modulo the $abc$-conjecture for number fields, that there exist infinitely many distinct monogenic cyclic polynomials of degree 4. A partial answer for the case $N=6$ is given in \cite{JonesEvenSextics}.

The second author observed that many of the monogenic cyclic polynomials appearing in the previously-cited articles occur as terms, or factors of terms, in certain ``classical" polynomial recurrence sequences \cite{Koshy,Horadam,Florez}, such as Fibonacci, Lucas, Vieta-Fibonacci, Vieta-Lucas etc.  However, in these sequences, some terms, and factors of terms, fail to be monogenic, cyclic or both. 
For example, in the Fibonacci polynomial sequence $\{F_n(x)\}_{n\ge 0}$, we see that $F_5(x)$ is monogenic but not cyclic, while the irreducible factor $x^6+7x^4+14x^2+7$ of $F_{14}(x)$ is cyclic but not monogenic. As a second example, in the Lucas sequence $\{L_n(x)\}_{n\ge 0}$, it can be shown that the irreducible terms are precisely the terms with index $2^k$, and these terms are both monogenic and cyclic of degree $2^k$. However, the irreducible factor $x^8+7x^6+14x^4+8x^2+1$ of $L_{15}(x)$ is neither monogenic nor cyclic. These observations provided the motivation for the current article. In particular, can a polynomial recurrence sequence be constructed such that the irreducible factors of all terms are monogenic and cyclic? In an effort to answer this question, we present a polynomial recurrence sequence in this article such that the irreducible factors of all terms are monogenic, and that there exist infinitely many terms for which all irreducible factors are also cyclic.

We require the following definition.
\begin{definition}\label{C}
  Let $N\ge 2$ be an integer, and let $\phi$ denote Euler's totient function. We say that {\em $N$ satisfies Condition $\CC$} if
  \[N\in \{2^a, \ p^a, \ 2p^a,\ 4p^a, \ p^aq^b, \ 2p^aq^b\},\]
   for some positive integers $a$ and $b$, and distinct odd primes $p$ and $q$ with $\gcd(\phi(p^a),\phi(q^b))=2$.
\end{definition}
The main result of this article is the following: 
\begin{theorem}\label{Thm:Main}\text{}
   Define the polynomial recurrence sequence $\{w_n(x)\}_{n\ge 0}$  as 
  \begin{align*}\label{Eq:w}
  \begin{split}
  w_0(x)&=1, \quad w_1(x)=1 \quad \mbox{and}\\
  w_n(x)&=(x-2)w_{n-1}(x)-w_{n-2}(x) \ \mbox{for $n\ge 2$.}
  \end{split}
  \end{align*}
  Then, for $n\ge 2$,
    \[w_{n}(x)=\ds \mathop{\prod_{d\mid (2n-1)}}_{d>1}\Omega_{d}(x),\] where  each $\Omega_d(x)$ is irreducible and monogenic of degree $\abs{\Gal(\Omega_d)}=\phi(d)/2$. Furthermore, $\Omega_d(x)$ is cyclic with $\Gal(\Omega_d)\simeq C_{\phi(d)/2}$ if and only if $d$ satisfies Condition $\CC$ as defined in Definition \eqref{C}. 
     \end{theorem}
   We then have the following corollary of Theorem \ref{Thm:Main}. 
   \begin{corollary}\label{Cor:Main}
   Let $n\ge 2$, and suppose that $n$ is such that $2n-1$ satisfies Condition $\CC$.
   Then,  every irreducible factor $\Omega_d(x)$ of $w_n(x)$ is both monogenic and cyclic. Moreover, $\Omega_{d_1}(x)$ and $\Omega_{d_2}(x)$ are distinct, where $d_i=2n_i-1$ with $d_1\ne d_2$ such that  each $n_i\ge 2$ and each $2n_i-1$ satisfies Condition $\CC$ as defined in Definition \eqref{C}.
   \end{corollary}

\section{Preliminaries}\label{Section:Prelim}
We begin with some basic information regarding the monogenicity of a polynomial. 
Suppose that $f(x)\in \Z[x]$ is monic and irreducible over $\Q$.
let $K=\Q(\theta)$ with ring of integers $\Z_K$, where $f(\theta)=0$. Then, the equation 
\begin{equation}\label{Eq:Dis-Dis}
\Delta(f)=\left[\Z_K:\Z[\theta]\right]^2\Delta(K)
\end{equation}
is well known  \cite{Cohen},
where $\Delta(f)$ and $\Delta(K)$ denote the discriminants over $\Q$, respectively, of $f(x)$ and the number field $K$. Recall that if $f(x)$ is monogenic, then $\Z_K=\Z[\theta]$. Thus, from \eqref{Eq:Dis-Dis}, we have that 
\begin{equation}\label{Eq:MonoDis}
f(x)\ \mbox{is monogenic if and only if} \ \Delta(f)=\Delta(K).
\end{equation} 
The following proposition will be useful in the proof of Corollary \ref{Cor:Main}.   
\begin{proposition}\label{Prop:equivalent}
  Let $f(x)$ and $g(x)$ be degree-$N$ monogenic cyclic polynomials. If $f(x)$ is equivalent to $g(x)$, then $\Delta(f)=\Delta(g)$.
\end{proposition}
\begin{proof} Suppose that $f(\alpha)=g(\beta)=0$. Let $K=\Q(\alpha)$ and $L=\Q(\beta)$. Since $f(x)$ and $g(x)$ are equivalent, we have that $K=L$, which implies that $\Delta(K)=\Delta(L)$. Since $f(x)$ and $g(x)$ are both monogenic, it follows from \eqref{Eq:MonoDis} that $\Delta(f)=\Delta(K)$ and $\Delta(g)=\Delta(L)$. Thus, $\Delta(f)=\Delta(g)$. 
 \end{proof}
 \begin{remark}
   The converse of Proposition \ref{Prop:equivalent} is false as can be seen from Example \ref{Ex}.
 \end{remark}

\begin{definition}\label{Def:Vieta}\cite{Horadam}
  We define  the \emph{Vieta-Lucas polynomials} $v_n(x)$ as 
  \[v_0(x)=2, \quad v_1(x)=x \quad \mbox{and} \quad v_n(x)=xv_{n-1}(x)-v_{n-2}(x) \quad \mbox{for all $n\ge 2$}.\]
  \end{definition}
\begin{remark}
  The Vieta-Lucas polynomials in Definition \ref{Def:Vieta} are also referred to simply as the \emph{Vieta polynomials} \cite{Robbins}.
\end{remark}

\begin{proposition}\label{Prop:Vieta}
 Let $v_n(x)$ denote the Vieta-Lucas polynomial of index $n\ge 1$. 
\begin{enumerate}
\item {\rm \cite{Horadam}} The zeros of $v_n(x)$ are
\[x=2\cos\left(\left(\dfrac{2j-1}{2n}\right)\pi\right), \quad j=1,2,\ldots, n.\]
Note if $2\nmid n$, then $x=0$ when $j=(n+1)/2$.
\item {\rm \cite{Robbins}}
\[v_n(x)=\sum_{j=0}^{\lfloor n/2\rfloor} (-1)^jB(n,j)x^{n-2j},\]
where
\[B(n,j)=\left(\frac{n}{n-j}\right)\binom{n-j}{j}=
\left\{\begin{array}{cl}
  1 & \mbox{if $j=0$}\\
  \left(\frac{n}{j}\right)\binom{n-j-1}{j-1} & \mbox{if $1\le j\le \lfloor n/2\rfloor$.}
\end{array}\right. \]
\item {\rm \cite{Florez}}\label{Vieta2} $\Delta(v_n)=2^{n-1}n^n$.
\end{enumerate}
\end{proposition}

\begin{corollary}\label{Cor:1}
  Let $n\ge 2$, and let $\zeta=e^{i\pi/(2(2n-1))}$. Let $W_n(x):=w_n(x^2)$.
  \begin{enumerate}
   \item \label{Cor:I1} $W_n(x)=v_{2n-1}(x)/x$;
   \item \label{Cor:I2} the $2n-2$ roots of $W_n(x)$ are
   \[RW=\left\{\zeta^{2j-1}+\zeta^{-(2j-1)}: j=1,2,\ldots, n-1,n+1,n+2,\ldots, 2n-1\right\};\]
  \item \label{Cor:I2A} the $n-1$ roots of $w_n(x)$ are
\[Rw=\left\{\zeta^{2(2j-1)}+\zeta^{-2(2j-1)}+2: j=1,2,\ldots, n-1\right\};\]
\item \label{Cor:I3} if $2n-1$ is prime, then $w_n(x)$ and $W_n(x)$ are both irreducible over $\Q$;
 \item \label{Cor:I4} $\Delta(w_n)=(2n-1)^{n-2}$ and $\Delta(W_n)=2^{2n-2}(2n-1)^{2n-3}$;
\end{enumerate}
 \end{corollary}
 \begin{proof}
 We use induction to prove item \eqref{Cor:I1}, which is readily verified for $n\in \{2,3\}$. Observe that $W_n(x)$ satisfies the recurrence
 \[W_0(x)=1,\quad W_1(x)=1 \quad \mbox{and} \quad W_n(x)=(x^2-2)W_{n-1}(x)-W_{n-2}(x) \quad \mbox{for $n\ge 2$.}\] Assume item \eqref{Cor:I1} is true for all $n\in \{2,\ldots ,k\}$ for some $k\ge 3$. Then, for $n=k+1$, we have
 \begin{align*}
   xW_{k+1}(x)&=x\left((x^2-2)W_{k}(x)-W_{k-1}(x)\right)\\
   &=(x^2-2)v_{2k-1}(x)-v_{2k-3}(x)\qquad \mbox{(by induction)}\\
   &=x^2v_{2k-1}(x)-2v_{2k-1}(x)-v_{2k-3}(x)\\
   &=x^2v_{2k-1}(x)-2(xv_{2k-2}(x)-v_{2k-3}(x))-v_{2k-3}(x)\\
   &=x\left(xv_{2k-1}(x)-v_{2k-2}(x)\right)-(xv_{2k-2}(x)-v_{2k-3}(x))\\
   &=xv_{2k}(x)-v_{2k-1}(x)\\
   &=v_{2k+1}(x).
    \end{align*}
    Item \eqref{Cor:I2} follows from item \eqref{Cor:I1} and Proposition \eqref{Prop:Vieta}. Since
    \[(\zeta^{2j-1}+\zeta^{-(2j-1)})^2=\zeta^{2(2j-1)}+\zeta^{-2(2j-1)}+2,\]
    item \eqref{Cor:I2A} follows from item \eqref{Cor:I1}, Proposition \eqref{Prop:Vieta} and the definition of $w_n(x)$. For Item \eqref{Cor:I3}, let $p=2n-1$ be prime. Then, for $1\le j\le (p-3)/2$, we see from Proposition \ref{Prop:Vieta} that the absolute value of the coefficient on $x^{(p-1)/2-j}$ in $w_n(x)$ is
  \[B(p,j)=\dfrac{p}{p-j}\binom{p-j}{j}\equiv 0 \pmod{p},\]
    and that the absolute value of $w_n(0)$ is
  \[B(p,(p-1)/2)=\left(\frac{p}{(p-1)/2}\right)\binom{(p-1)/2}{(p-3)/2}=p.\]
  We conclude that both $w_n(x)$ and $W_n(x)$ are $p$-Eisenstein, and consequently irreducible over $\Q$. Item \eqref{Cor:I4} follows from Proposition \ref{Prop:Vieta}, \cite[Theorem 2.7]{HJMCC} and \cite{Janson}.
  \end{proof}

\section{The Proof of Theorem \ref{Thm:Main}}\label{Section:MainProof}
We first prove two lemmas. The first lemma, which does not seem to appear in the literature, is of some interest in its own right.
\begin{lemma}\label{Lem:Max}
  Let $N\ge 3$ be an integer, and let $\zeta$ be a primitive $N$th root of unity. 
  Then the Galois group of the maximal real subfield $\Q(\zeta+\zeta^{-1})$ of $\Q(\zeta)$ is cyclic if and only if $N$ satisfies Condition $\CC$ as defined in Definition \eqref{C}.
  \end{lemma}
\begin{proof}
We give details only for the case when $2\nmid N$ since the case $2\mid N$ is similar. Let $N=\prod_{i=1}^kp_i^{e_i}$ be the factorization of $N$ into distinct odd prime powers. Then,
\begin{align}\label{G}
\nonumber G:=\Gal(\Q(\zeta))&\simeq (C_N)^*\\
\nonumber &\simeq C_{p^{e_1}}^*\times \cdots \times C_{p^{e_k}}^*\\
& \simeq C_{p^{e_1-1}(p_1-1)}\times \cdots \times C_{p^{e_k-1}(p_k-1)}.
\end{align}
Since the fixed field of the complex conjugation automorphism $\sigma:\zeta \mapsto \zeta^{-1}$ of order 2 in $G$ is $\Q(\zeta+\zeta^{-1})$, it follows that
$\sigma$ embeds into $\Gal(\Q(\zeta))$ in \eqref{G} as $(\sigma_1,\sigma_2,\ldots, \sigma_k)$, where $\sigma_i$ is the unique element of order 2 in $C_{p^{e_i-1}(p_i-1)}$. Thus, with $H=\left<\sigma\right>\simeq \left<(\sigma_1,\sigma_2,\ldots, \sigma_k)\right>$, we have that
  \begin{equation}\label{G/H}
  G/H\simeq \Gal(\Q(\zeta+\zeta^{-1})). 
  \end{equation}
    It is then easy to verify that
      \[\underbrace{(\sigma_1,1,\ldots, 1)}_\text{$k$-tuple}\left<\sigma\right> \quad \mbox{and} \quad \underbrace{(1,\sigma_2,1,\ldots, 1)}_\text{$k$-tuple}\left<\sigma\right>\]
   are distinct elements of order 2 in $G/H$ if $k\ge 3$, so that $G/H$ cannot be cyclic in this situation. Hence, if $G/H$ is cyclic, then  $k\le 2$.
      If $k=1$, then $G/H$ is cyclic since $G$ is cyclic.

   Thus, we only need to examine the case $k=2$:
   \[G\simeq A_1\times A_2 \quad \mbox{and} \quad  G/H\simeq \left(A_1\times A_2\right)/H,\] where
   $A_i=C_{p_i^{e_i-1}(p_i-1)}$ and $H=\left<(\sigma_1,\sigma_2)\right>$. 
   Let $\Gamma:=\gcd(\phi(p_1^{e_1}),\phi(p_2^{e_2}))$ and let $2^{r_i}\mid\mid (p_i-1)$. Suppose that $\min_{i}(r_i)=2$, so that $2^2\mid \Gamma$. Then there exists $(\gamma_1,\gamma_2)\in G$ with $\ord_{A_i}(\gamma_i)=4$, where $\ord_{A_i}(\gamma_i)$ denotes the order of $\gamma_i$ in $A_i$. Straightforward computations reveal
     \begin{gather*}
     (\gamma_1^2,\gamma_2^2)H=(\sigma_1,\sigma_2)H=H, \quad (\sigma_1,1)H      \ne H, \quad (\gamma_1,\gamma_2)H \ne (\sigma_1,1)H\\ 
     \mbox{and} \ \ord_{G/H}\left((\sigma_1,1)H\right)     =\ord_{G/H}\left((\gamma_1,\gamma_2)H\right)=2,
     \end{gather*}
      which implies that $G/H$ contains two distinct elements of order 2, namely $(\sigma_1,1)H$ and $(\gamma_1,\gamma_2)H$, and therefore $G/H$ is not cyclic. Hence, $\min_{i}(r_i)=1$ if $G/H$ is cyclic, and in this case, we have that $2\mid \mid \Gamma$.

        Observe that if $e_1>1$ and $p_1\mid (p_2-1)$, then $G/H$ contains a subgroup isomorphic to $C_{p_1}\times C_{p_1}$, so that $G/H$ is not cyclic. In this case we have that $p_1\mid \Gamma$. Moreover, if $e_1=1$ or $p_1\nmid (p_2-1)$, then the $p_1$-subgroup of $G/H$ is cyclic and $p_1\nmid \Gamma$. Reversing the roles of $p_1$ and $p_2$ yields similar results for $p_2$. Combining all of this information completes the proof of the lemma.
\end{proof}

The next lemma calculates the discriminant of the maximal real subfield of a cyclotomic field.
\begin{lemma}\label{Lem:DiscMax}
  Let $N\ge 3$ be an integer. Let $\zeta$ denote a primitive $N$th root of unity, and let $K=\Q(\zeta+\zeta^{-1})$. Then
  \[\Delta(K)=\left\{\begin{array}{cl}
    p^{(p^{k-1}pk-k-1)-1)/2} & \mbox{if $N\in \{p^k,2p^k\}$, where $p\ge 3$ is prime;}\\[3pt]
    2^{2^{k-2}(k-1)-1} & \mbox{if $N=2^k$, where $k\ge 2$;}\\[3pt]
    \dfrac{N^{\phi(N)/2}}{\ds \prod_{p\mid N}p^{\phi(N)/(2(p-1))}} & \mbox{otherwise.}
  \end{array}\right.\]
\end{lemma}
\begin{proof}
  Let $f(x)=x^2-(\zeta+\zeta^{-1})x+1$, the minimal polynomial for $\zeta$ over $K$. Let $\NN:=\NN_{L/\Q}$ denote the norm, where $L=\Q(\zeta)$. Then
  \begin{align*}
    \abs{\Delta(L)}&=\Delta(K)^2\abs{\NN(f'(\zeta))} \quad \mbox{(from \cite{Neukirch})}\\
    &=\Delta(K)^2\abs{\NN(\zeta-\zeta^{-1})}\\
    &=\Delta(K)^2\abs{\NN(\zeta-1)\NN(\zeta+1)\NN(\zeta)^{-1}}\\
    &=\left\{\begin{array}{cl}
      p\Delta(K)^2 & \mbox{if $N\in \{p^k,2p^k\}$, where $p\ge 3$ is prime;}\\[3pt]
      4\Delta(K)^2 & \mbox{if $N=2^k$;}\\[3pt]
      \Delta(K)^2 & \mbox{otherwise.}
    \end{array}\right.
  \end{align*} The lemma then follows from the formula for $\Delta(L)$ \cite[Proposition 2.7]{Washington}, and the fact that $\Delta(K)>0$ since $K$ is a totally real number field.
\end{proof}

\begin{proof}[Proof of Theorem \ref{Thm:Main}]\text{}
Let $\zeta=e^{i\pi/(2(2n-1))}$.
For each divisor $d>1$ of $2n-1$ and $j\in \{1,2,\ldots, n-1\}$, define
\[Rw_d:=\left\{\zeta^{2(2j-1)}+\zeta^{-2(2j-1)}+2:\ \gcd(2j-1,2n-1)=(2n-1)/d\right\}.\] Since $\gcd(2j-1,2n-1)=(2n-1)/d$, we have that $m=(2j-1)d/(2n-1)$ is an odd integer and $\gcd(m,d)=1$. Thus,
$\abs{Rw_d}=\phi(d)/2$, and $Rw_a \cap Rw_b=\varnothing$ if $a\ne b$. Furthermore, from Corollary \ref{Cor:1}, it follows that
$Rw=\bigcup_{d}Rw_d$.  Define
\[\Omega_d(x):=\prod_{r_i\in Rw_d}(x-r_i).\] Observe that $\rho:=\zeta^{2(2n-1)/d}$ is a primitive $2d$th  root of unity. Then, $\Omega_d(x)\in \Z[x]$ is irreducible over $\Q$ of degree $\phi(2d)/2=\phi(d)/2$, and the splitting field of $\Omega_d(x)$ is
\[\Q(\rho+\rho^{-1}+2)=\Q(\rho+\rho^{-1}),\] which is the maximal real subfield of $\Q(\rho)$. Hence, $\Omega_d(x)$ is monogenic by \cite[Proposition 2.16]{Washington}.
By Lemma \ref{Lem:Max}, $\Omega_d(x)$ is cyclic with $\Gal(\Omega_d)\simeq C_{\phi(d)/2}$ if and only if $d$ satisfies Condition $\CC$.
\end{proof}

\section{The Proof of Corollary \ref{Cor:Main}}\label{Section:MainCorProof}
\begin{proof}
  If $n\ge 2$ is an integer such that $2n-1$ satisfies Condition $\CC$, then every divisor $d$ of $2n-1$ satisfies Condition $\CC$. Hence, every irreducible factor $\Omega_d(x)$ of $w_n(x)$ is both monogenic and cyclic by Theorem \ref{Thm:Main}.

  For the second part, suppose that $n_1\ge 2$ and $n_2\ge 2$ are integers such that both $d_1:=2n_1-1$ and $d_2:=2n_2-1$ satisfy Condition $\CC$. Then, by Theorem \ref{Thm:Main}, both $\Omega_{d_1}(x)$ and $\Omega_{d_2}(x)$ are monogenic and cyclic. For $i\in \{1,2\}$, let $K_i=\Q(\alpha_i)$, where $\Omega_{d_i}(\alpha_i)=0$. Suppose that $\Omega_{d_1}(x)$ and $\Omega_{d_2}(x)$ are equivalent. Then, noting Proposition \ref{Prop:equivalent}, we have that
  \begin{gather*}
  \phi(d_1)=\phi(d_2), \quad K_1=K_2, \quad \Delta(\Omega_{d_i})=\Delta(K_i)\\
   \mbox{and} \ \Gal(\Omega_{d_1})\simeq\Gal(\Omega_{d_2})\simeq C_{\phi(d_i)/2}.
   \end{gather*} Consequently, $n_1=n_2$ by  Lemma \ref{Lem:DiscMax}.
  \end{proof}

\section{Final Comments}\label{Section:Final Comments}
For any $n\ge 2$, we see from Theorem \ref{Thm:Main} that $\Omega_{2n-1}(x)$ is a monogenic divisor of $w_n(x)$ such that $\Omega_{2n-1}(x)$ does not divide $w_m(x)$ for any $m<n$. We call $\Omega_{2n-1}(x)$ a {\em primitive divisor} of $w_n(x)$. In fact, $\Omega_{2n-1}(x)$ is the unique primitive divisor of $w_n(x)$ since every other divisor $\Omega_d(x)$ of $w_n(x)$, where $1<d<2n-1$ is a divisor of $2n-1$, must divide $w_{(d+1)/2}(x)$.  Moreover, every maximal real subfield $K_{N}^{+}$ of a cyclotomic field $K_{N}$ of odd index $N$ is realized as the splitting field of the primitive divisor $\Omega_{2n-1}(x)$ of $w_n(x)$ for some $n$. In other words, for a given odd integer $N\ge 3$, the primitive divisor $\Omega_{2n-1}(x)=\Omega_{N}(x)$ of $w_n(x)$ is a monogenic polynomial that generates $K_N^{+}$; and is also cyclic, provided $N$ satisfies Condition $\CC$.

\subsection*{Acknowledgements}
The authors thank the anonymous referee for the suggestions that helped to improve the paper. 


\end{document}